\documentclass[11pt]{amsart}
\usepackage[english]{babel}
\usepackage[latin1]{inputenc}
\usepackage[dvips,final]{graphics}
\usepackage{amsmath,amsfonts,amssymb,amsthm,amscd,array,stmaryrd,mathrsfs}
\usepackage{url}
 \usepackage[all]{xy}
\usepackage{graphicx}
\usepackage{color}              
\usepackage{epsfig}
\usepackage{psfrag}
\usepackage{epstopdf}

\usepackage{textcomp}

\theoremstyle{plain}
\newtheorem{thm}{Theorem}
\newtheorem{lem}{Lemma}[section]

\theoremstyle{definition}

\newtheorem{ex}[lem]{Example}

\newcommand{\R}{\mathbb{R}}

\newcommand{\C}{\mathbb{C}}
\newcommand{\RP}{\mathbb{RP}}
\newcommand{\CP}{\mathbb{CP}}
\newcommand{\HP}{\mathbb{HP}}
\newcommand{\OP}{\mathbb{OP}}
\newcommand{\bbH}{\mathbb{H}}

\newcommand{\bbO}{\mathbb{O}}

\newcommand{\N}{\mathbb{N}}

\def\r{\rho}

\begin{document}

\title[Affine Hopf fibrations]{Affine Hopf fibrations}

\author{
Valentin Ovsienko
\and
Serge Tabachnikov}

\address{
Valentin Ovsienko,
CNRS,
Laboratoire de Math\'ematiques, 
Universit\'e de Reims-Champagne-Ardenne, 
FR 3399 CNRS, F-51687, Reims, France}
\email{valentin.ovsienko@univ-reims.fr}

\address{
Serge Tabachnikov,
Department of Mathematics,
Pennsylvania State University, 
University Park, PA 16802, USA
}
\email{tabachni@math.psu.edu}

\date{}


\maketitle

\thispagestyle{empty}

\section{Introduction}

There exist four fiber bundles, called the {\it Hopf fibrations}, whose
fibers, total spaces, and bases are spheres:
\begin{equation}
\label{CHFib}
 \xymatrix{
S^0\ar@{->}[r]&S^1\ar@<-1pt>@{->}[d]\\
&S^1
}
\qquad
 \xymatrix{
S^1\ar@{->}[r]&S^3\ar@<-1pt>@{->}[d]\\
&S^2
}
\qquad
 \xymatrix{
S^3\ar@{->}[r]&S^7\ar@<-1pt>@{->}[d]\\
&S^4
}
\qquad
 \xymatrix{
S^7\ar@{->}[r]&S^{15}\ar@<-3pt>@{->}[d]\\
&\!\!S^8.
}
\end{equation}
The vertical arrows denote projections whose fibers,
i.e., inverse images of points,
are represented by the horizontal arrows.
The term ``fibration'' means that, locally, the total space is the product of the base and the fiber, hence the bigger spheres are filed with smaller ones.
For instance, through every point of $S^3$ there passes
one circle, $S^1$, and different circles do not intersect.

The construction of Hopf fibrations is very simple (we recall it below)
and is based on
the algebras $\R,\C,\bbH,\bbO$; according to a famous theorem of Adams~\cite{Ada1},
that there are no such fibrations in other dimensions.
Note that the ``magic" numbers $1,2,4,8$ appear twice: 
as the dimensions of the bases, and also, this time shifted by~$1$,
as the dimensions of the fibers.

In a sense, contemporary algebraic topology has grown up with the Hopf fibrations:
the development of   the theory of
characteristic classes, homotopy theory,  and  $K$-theory
was much influenced by the study of Hopf fibrations;
see \cite{Ati,Bot,num,Hus,Kar}.

This note, based on a recent work~\cite{OT}, concerns a similar problem: for
which pairs $(p,n)$, is $\R^n$ foliated by {\it pairwise skew} affine subspaces $\R^p$?
Two disjoint affine subspaces are pairwise skew if they do not contain parallel lines.\footnote{Note that if $n$-dimensional space admits pairwise skew $p$-dimensional subspaces then $n\ge 2p+1$.}
In other words, we consider  fiber bundles
\begin{equation}
\label{AffHFib}
 \xymatrix{
\R^p\ar@{->}[r]&\R^n\ar@<-1pt>@{->}[d]\\
&\;\;\;\R^{n-p}
}
\end{equation}
with pairwise skew affine fibers. We call them  \textit{affine Hopf fibrations}.
We shall give an answer in terms of the Hurwitz-Radon function 
that appeared in topology in the theorem 
of Adams~\cite{Ada} about vector fields on spheres.

The reason to consider such ``skew" foliations is two-fold. On the one hand, the central projection of a Hopf bundle from a sphere to Euclidean space yields such a foliation. And, on the other hand, without the skew restriction, the problem is trivial: for any $p<n$, one can foliate $\R^n$ by parallel copies of $\R^p$.

\section{The Hurwitz-Radon function}

The sequence
$$
\mathbf{1,2,4,8},9,10,12,16,17,18,20,24,25,26,28,32\ldots
$$
is known (A003485, in the Sloane encyclopedia of integral sequences \cite{Slo}) as
the {\it Hurwitz-Radon function evaluated at powers of $2$}.
Namely, the Hurwitz-Radon function $\rho:\N\to\N$
is defined as follows.
Every natural number $N$ can be written as  $N=2^n(2m+1)$, and the
function $\rho$ depends only on the dyadic part of $N$, that is,
$$
\r(N)=\r(2^n),
$$
and
$$
\r(N)=\left\{ \begin{array}{lcll}
2n+1, \quad&n\equiv& 0 &\mod 4\\
2n, \quad&n\equiv& 1,2 &\mod 4\\
2n+2, & n\equiv &3 &\mod 4.
\end{array}
\right.
$$

The numbers $N=1,2,4,8$ are the only numbers for which $\r(N)=N$.
The above  formula can be equivalently written as
$$
\r(2^{n+4})=\r(2^n)+8,
$$
which is easier to remember.

\subsection{The history of the Hurwitz-Radon function}
The Hurwitz-Radon function was discovered around 1920, independently, by
Adolf Hurwitz and Jean~Radon~\cite{Hur,Rad}\footnote{The paper of Hurwitz was published posthumously.}.
Both were working on the problem of ``square identities'',
or composition of quadratic forms, that is, formulas of the type
$$
(a^2_1+\cdots+a^2_r) (b^2_1+\cdots+b^2_s)=c^2_1+\cdots+c^2_N,
$$
where $c_1,\ldots,c_N$ are bilinear forms in $a_1,\ldots,a_r$ and $b_1,\ldots,b_s$
with real coefficients.
The above identity is said to be of size $[r,s,N]$.

For example, a formula of size $[4,4,4]$ is Euler's four-square identity
\begin{equation*}
\begin{array}{rcl}
(a_1^2+a_2^2+a_3^2+a_4^2)(b_1^2+b_2^2+b_3^2+b_4^2)&=&
(a_1 b_1 + a_2 b_2 + a_3 b_3 + a_4 b_4)^2
\\[2pt]
&+&
(a_1 b_2 - a_2 b_1 + a_3 b_4 - a_4 b_3)^2 
\\[2pt]
&+&
(a_1 b_3 - a_2 b_4 - a_3 b_1 + a_4 b_2)^2
\\[2pt]
&+&
(a_1 b_4 + a_2 b_3 - a_3 b_2 - a_4 b_1)^2.
\end{array}
\end{equation*}
that corresponds to multiplication of quaternions.
Euler's identity is nothing other than the property 
$\vert a\vert\vert b\vert=\vert ab\vert$ in $\bbH$.
It was discovered by Leonhard Euler in~1748,
almost a hundred years before the discovery of quaternions.
Note also that a similar formula of size $[2,2,2]$, that corresponds
to multiplication of complex numbers, was known to Diophantus,
and also appeared in the early VII century in a book of
an Indian mathematician Brahmagupta;
it was also used by Fibonacci in his ``Book of squares''.
There exists also a square identity of size $[8,8,8]$, found by
a Dutch mathematician Ferdinand Degen in 1818, that corresponds
to multiplication of octonions.

In 1898, Hurwitz proved his famous $1,2,4,8$ theorem,
stating that  a square identity of size $[N,N,N]$ can exist only for 
$n\in\{1,2,4,8\}$.
This is equivalent to the statement that $\R,\C,\bbH,\bbO$ are
the only real normed division algebras with unit.
He also formulated a problem to characterize triplets
$r,s,N\in\N$ for which there exists a square identity of size $[r,s,N]$.
This problem remains widely open.
Traditionally considered as a problem of number theory, it plays important role
in many other areas of mathematics, for more details see~\cite{Sha}.

Hurwitz and Radon proved that a formula of size $[r,N,N]$ exists if and only if $r\leq \r(N)$,
and this is still the only case where the Hurwitz problem is solved. 

The Hurwitz-Radon function later appeared, sometimes unexpectedly, 
in many different areas, such as 
algebra, representation theory, geometry, topology, and combinatorics.
The function gives the dimension of irreducible representations
of Clifford algebras that are important in mathematics and physics.
Recently the Hurwitz-Radon function appeared in the context
of multi-antennas wireless communication.

The ``exceptional'' numbers $1,2,4,8$ (or $0,1,3,7$)
often appear in mathematics and physics.
The best known subjects are:
the celebrated results of topology about parallelizable spheres and Hopf fibrations,
and the famous theorem of algebra about classification of real normed division algebras. 
In these subjects, the above numbers appear as the critical dimensions.

\subsection{The Adams theorem about vector fields on spheres}

The theorem of Adams~\cite{Ada} is perhaps the best known, 
and one of the most beautiful applications
of the Hurwitz-Radon function.

\begin{thm}
\label{Athm}
 The maximal number of vector fields on the sphere $S^{N-1}$,
linearly independent at each point, is $\r(N)-1$.
\end{thm}

The existence part of this theorem follows from the Hurwitz-Radon construction of square identities. 
Suppose that we have a square identity of size $[r,N,N]$.
Consider the $N$-vector $\bar c=(c_1,\ldots,c_N)^T$.
Since every $c_i$ is a bilinear form in $a$'s and $b$'s,
we have:
\begin{equation}
\label{MatHR}
\bar c=
\left(a_1A_1+\cdots+a_rA_r\right)
\left(
\begin{array}{c}
b_1\\
\vdots\\
b_N
\end{array}
\right),
\end{equation}
where $A_1,\ldots,A_r$ are $N\times{}N$-matrices with real coefficients.
The existence of a square identity implies that, for every non-zero $\bar b\in\R^N$,
the vectors $A_1\bar b,\ldots,A_r\bar b$ are linearly independent,
i.e., that $\bar c\not=0$ in~(\ref{MatHR}) whenever $\bar a\in\R^r$ is non-zero.
Indeed, the right-hand side of the square identity is the square of the norm of
$\bar c$, while the left-hand side is the product of norms of $\bar a$ and~$\bar b$.

Therefore, the existence of a square identity of size $[r,N,N]$ guarantees the existence
of $r$ (linear) vector fields in $\R^N\setminus\{0\}$, linearly independent at each point.
To obtain (at least) $r-1$ independent vector fields on $S^{N-1}$,
one now restricts the constructed $r$ fields to the round sphere and projects
to its tangent plane.

Let us mention that the existence part of the Adams theorem was
 known since 1920's;
the converse statement is one of the most difficult results of algebraic  topology.

Adams' theorem generalizes the theorem that 
the spheres $S^0,S^1,S^3$ and $S^7$ are the only {\it parallelizable}\footnote{
Parallelizability of a manifold means that it admits as many vector fields linearly independent
at each point as its dimension. 
Examples of parallelizable manifolds are Lie groups.} spheres.
This was proved independently in 1958 by Kervaire~\cite{Ker} and by Bott and Milnor~\cite{BM}.

\section{The classical Hopf fibrations}\label{HS3}

\subsection{The Hopf fibration of $S^3$}\label{H3CF}

The  fibration
$$
 \xymatrix{
S^1\ar@{->}[r]&S^3\ar@{->}[d]\\
&S^2
}
$$
was defined and studied by Heintz Hopf in 1931, see \cite{Hop}.

The fibers of the Hopf fibration
are obtained by the action of $S^1$, the   group of complex numbers with unit absolute values, in 
the complex plane $\C^{2}$.
The action is given by the formula
$$
e^{i\phi}\left(z_1,z_{2}\right)=\left(e^{i\phi}z_1,e^{i\phi}z_{2}\right).
$$
Restricting this action to the unit sphere $S^3\subset\C^{2}$
given by the equation 
$$
|z_1|^2+|z_2|^2=1,
$$
one obtains a fibration of $S^3$ by circles.

Equivalently, one may consider the projection
$\C^2 - \{O\}\to\CP^1$ that sends a point of the complex plane
to the line through this point and the origin.
In coordinates,  
$$
\left(z_1,z_{2}\right)\mapsto\left(z_1:z_{2}\right).
$$
Restricting this projection to $S^3\subset\C^{2}$, and noticing that
$\CP^1\simeq{}S^2$, the Riemann sphere, one obtains the Hopf fibration.

\begin{figure}[hbtp]
\vskip-8cm
\includegraphics[width=13cm]{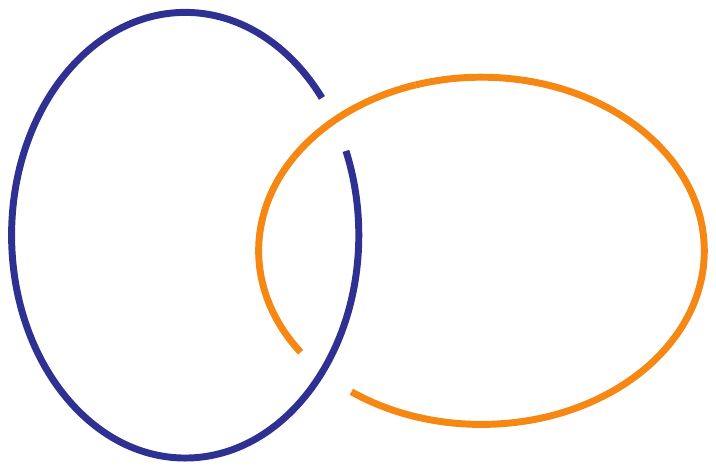}
\vskip-7cm
\caption{The Hopf link: two fibers of the Hopf fibration.}
\label{LinkC}
\end{figure}

Let us mention three properties of the Hopf fibration:
\begin{enumerate}
\item every two fibers of the Hopf fibration are linked non-trivially, see Figure \ref{LinkC};
\item every fiber is a {great circle} in $S^3$;
\item these great circles are equidistant from each other.
\end{enumerate}

For visualization of the Hopf fibration, we recommend the movies \cite{Dim,Vid}.

\subsection{Higher dimensional Hopf fibrations}

The   construction of the Hopf fibration of $S^3$  has two 
immediate generalizations.

Let us replace the algebra of complex numbers $\C$ 
by the other division algebras, $\R,\bbH$ and~$\bbO$.
The projections from the quaternionic and octonionic
planes to the respective projective lines, $\bbH^2\to\HP^1$ and $\bbO^2\to\OP^1$,
are defined as in the complex case.
Restricted to the spheres $S^7\subset\bbH^2$ and $S^{15}\subset\bbO^2$,
this projections give rise to the following fibrations
$$
 \xymatrix{
S^7\ar@{->}[r]&\HP^1
},
\qquad
 \xymatrix{
S^{15} \ar@{->}[r]&\OP^1
 }
$$
whose fibers are  $S^3$ and $S^7$, respectively
(the quaternions and octonions of absolute value $1$).
Note that  $\HP^1\simeq{}S^4$ and $\OP^1\simeq{}S^8$,
since these are quaternionic and octonionic lines completed by one point at infinity.\footnote{See \cite{CS} concerning octonions and $\OP^1$.}
Finally, replacing $\C$ by $\R$, one obtains the fibration
$
S^0\to{}S^1\to\RP^1
$
whose fibers consist of two points. This provides a complete list of Hopf fibrations (\ref{CHFib}).

Another straightforward generalization of the Hopf fibrations
is obtained replacing the plane by  $(m+1)$-dimensional complex vector space.
Consider the action of $S^1$ in  $\C^{m+1}$ given by
$$
e^{i\phi}\left(z_1,z_{2},\ldots,z_{m+1}\right)=
\left(e^{i\phi}z_1,e^{i\phi}z_{2},\ldots,e^{i\phi}z_{m+1}\right).
$$
Restricting this action to the sphere $S^{2m+1}\subset\C^{m+1}$,
one obtains a fibration of $S^{2m+1}$ by $S^1$. The base of this fibration is the complex projective space $\CP^m$.

The same construction can be applied in the quaternionic cases,
replacing $S^1$ by the unit sphere $S^3\subset\bbH$.
One  obtains three infinite series of  fibrations:
\begin{equation}
\label{HHFib}
 \xymatrix{
S^0\ar@{->}[r]&S^n\ar@<-2pt>@{->}[d]\\
&\RP^n
}
\qquad
 \xymatrix{
S^1\ar@{->}[r]&S^{2m+1}\ar@<-3pt>@{->}[d]\\
&\CP^n
}
\qquad
 \xymatrix{
S^3\ar@{->}[r]&S^{4m+3}\ar@<-3pt>@{->}[d]\\
&\HP^n
}
\end{equation}
This construction
does not work in the octonionic case, and  the fibration with fibers $S^7\to{}S^{15}\to{}S^8$ remains exceptional, and does not belong to any series.

It is still true that every fiber of a fibration~(\ref{HHFib}) is a great sphere.

\subsection{Projection to Euclidean spaces}

\begin{figure}[hbtp]
\includegraphics[width=12cm]{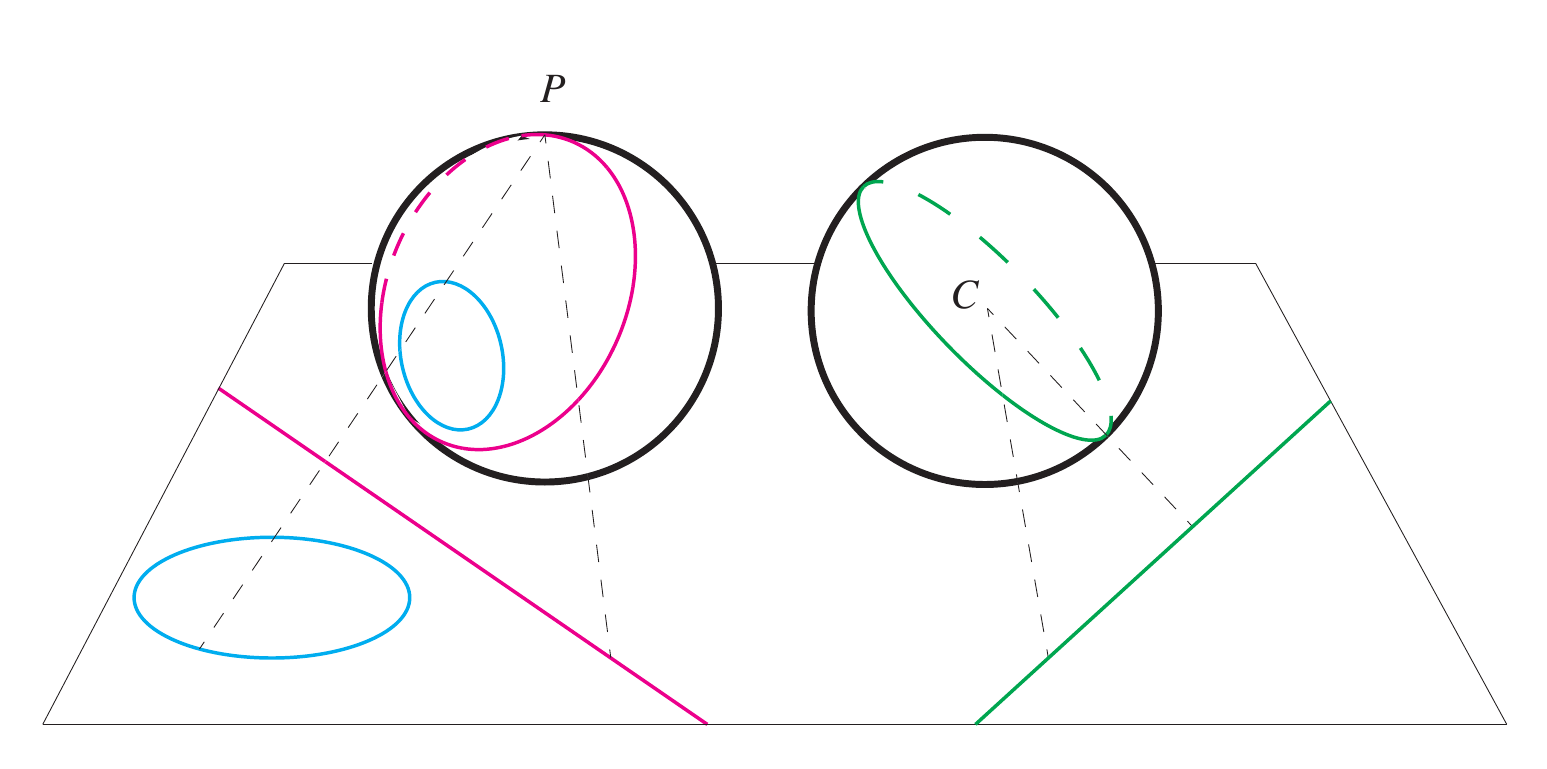}
\caption{The stereographic projection and the central projection: 
conformal geometry {vs} projective geometry.}
\label{ProjF}
\end{figure}

The usual way to ``visualize'' the Hopf fibration is by way of
projecting them to the Euclidean spaces.

The {\it stereographic projection}
is ubiquitous. This map identifies the sphere with one point, the North Pole $P$, deleted, 
with $n$-dimensional Euclidean space (the tangent space at the South pole), see Figure \ref{ProjF}.
The main property of the stereographic projection is that
this is a conformal map: it preserves the angles.
It sends circles to circles, except for the circles that pass through~$P$:  
they are sent to straight lines.

An equally important  projection from the $n$-dimensional sphere to the
$n$-dimensional Euclidean space is the {\it central projection}.
It projects the sphere $S^n \subset \R^{n+1}$ with the removed equatorial sphere $S^{n-1}$
from the center of the sphere~$C$ to an affine space $\R^n \subset \R^{n+1}$ that does not pass through the origin. The central projection sends  great circles to straight lines,
and  great spheres to affine planes.

\begin{figure}[hbtp]
\includegraphics[width=4cm]{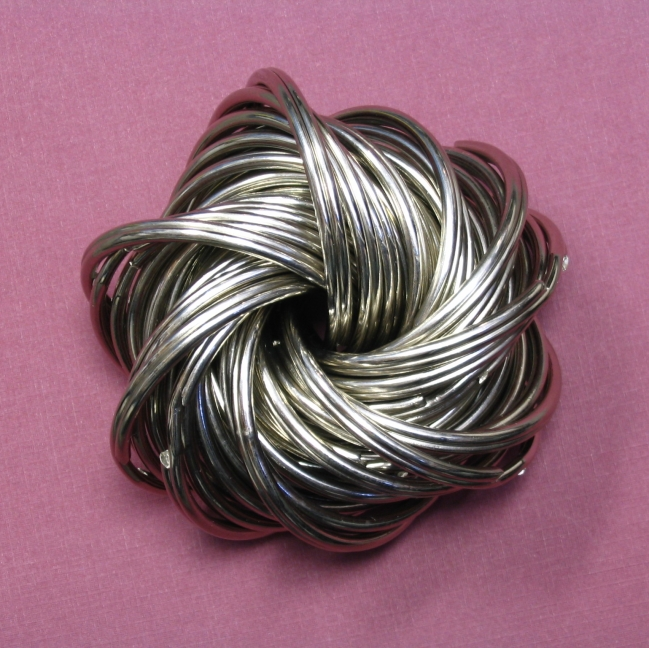}
\caption{A key-ring model of Villarceau circles, Wikimedia Commons.}
\label{keyrings}
\end{figure}

\begin{figure}[hbtp]
\includegraphics[width=10cm]{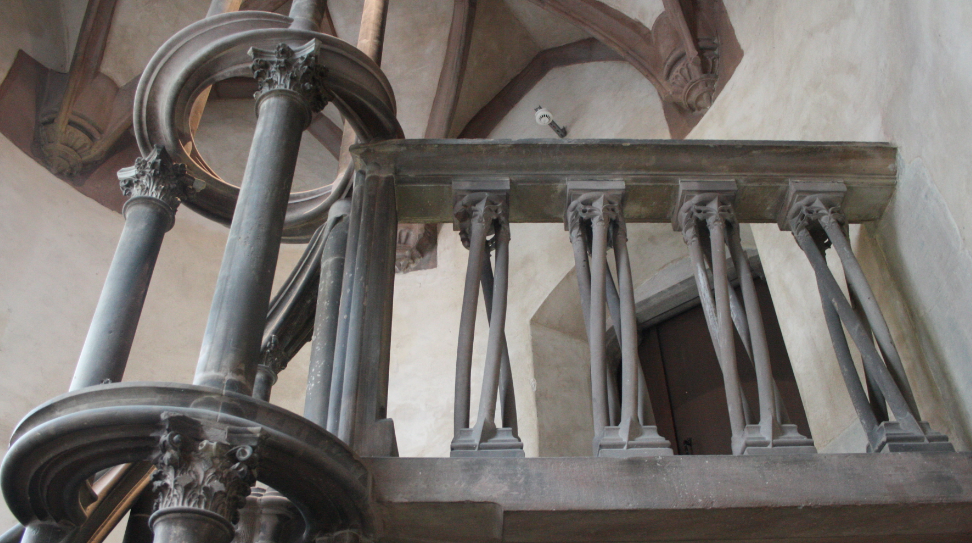}
\caption{Villarceau circles and skew lines 
at Mus\'ee de l'Oeuvre Notre-Dame, Strasbourg, a photograph by Luba Shenderova-Fock.}
\label{Villa}
\end{figure}

The image of the Hopf fibration of $S^3$ under the stereographic projection is a beautiful 
geometric structure in $\R^3$ called the Villarceau circles. These circles foliate $\R^3$ with a straight line removed (this line is the stereographic image of the Hopf circle through point $P$). 
See Figures~\ref{keyrings} and~\ref{Villa}.

\begin{figure}[hbtp]
\includegraphics[width=5cm]{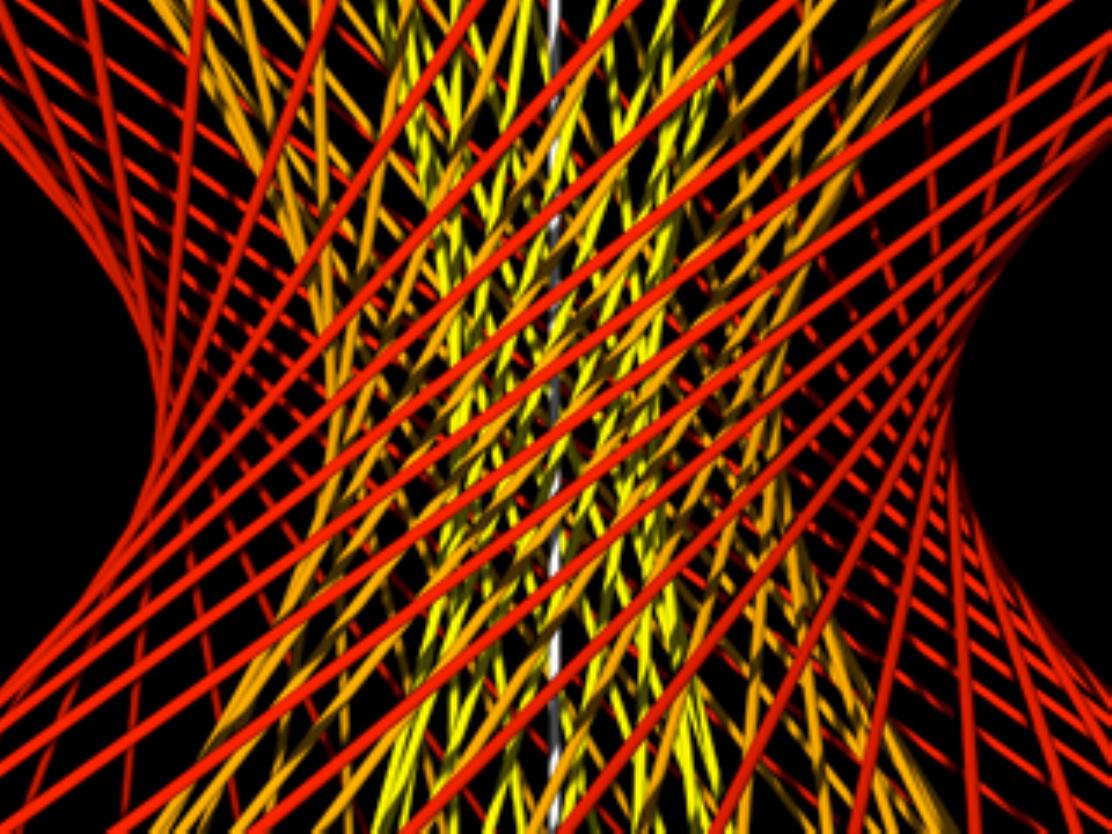}
\caption{``Nested hyperboloids" by David Eppstein, from the Wikipedia article
``Skew lines''.}
\label{Hof}
\end{figure}

The image of the Hopf fibration of $S^3$ under the central projection is a fibration of $\R^3$ by lines, see Figure \ref{Hof}.
 These lines are pairwise skew: two parallel lines `intersect' at infinity, and their preimages in $S^3$ are great circles that intersect at a point on the equator. 

\begin{ex}
\label{PHFi}
The central projections of the fibrations~(\ref{CHFib}) and~(\ref{HHFib}) 
give fibrations of the Euclidean space by pairwise skew affine subspaces:
$$
 \xymatrix{
\R^p\ar@{->}[r]&\R^n\ar@<-1pt>@{->}[r]
&\R^{n-p},
}
\qquad
\hbox{where}
\qquad
(p,n)=\left\{
\begin{array}{l}
(1,\,2m+1),\\
(3,\,4m+3),\\
(7,\,15).
\end{array}
\right.
$$
\end{ex}

\section{Affine Hopf fibrations}

Let us consider the affine Hopf fibrations~(\ref{AffHFib})
in full generality.

\subsection{Existence theorem and first examples}\label{ExSec}

The problem of describing all pairs $(p,n)$ for which such fibrations exist
was  solved in~\cite{OT}.\footnote{The next step would be a classification of affine Hopf fibrations. 
See \cite{Sal,Har} for partial results. 
A classification of fibrations of $S^3$ by great circles was obtained by 
H. Gluck and F. Warner \cite{GW}: 
they are in one-one correspondence with distance decreasing self-maps of a 2-sphere.}

\begin{thm}
\label{TheThm}
$\R^n$ admits a  fibration by $p$-dimensional pairwise skew affine subspaces
if and only if 
$
p\leq\r(n-p)-1.
$
\end{thm}

The condition of this theorem is satisfied in Examples \ref{PHFi}.
Indeed, if $p=1$ and $n=2m+1$, then $n-p$ is even, so that $\rho(n-p)\geq2$.
If $p=3$ and $n=4m+3$, then $n-p$ is 
a multiple of $4$, so that $\rho(n-p)\geq4$.
Finally, if $p=7$ and $n=15$, then $\rho(n-p)=8$.

Let us sketch the proof of Theorem \ref{TheThm}.
The existence part is, again, a consequence of the Hurwitz-Radon construction.
Given a square identity of size $[r,N,N]$, rewrite the vector of bilinear forms
$\bar c=(c_1,\ldots,c_N)^T$ of the right-hand-side in the following form
$$
\bar c=
\left(b_1B_1+\cdots+b_NB_N\right)
\left(
\begin{array}{c}
a_1\\
\vdots\\
a_r
\end{array}
\right),
$$
where $B_1,\ldots,B_N$ are $N\times{}r$-matrices with real coefficients.
Using the notation 
$$
B(\bar b):=b_1B_1+\cdots+b_NB_N,
$$
the square identity implies that
$B(\bar b)$ is of {\it maximal rank} (i.e., of rank $r$) for all $\bar b\not=0$.
This matrix expression is of course equivalent (actually, dual) to~(\ref{MatHR}).

Furthermore, using a linear change of coordinates $(b_1,\ldots,b_N)$, we can assume that
the last column of each matrix $B_i$ is the vector
$(0,\ldots,0,1,0,\ldots,0)^T$
with $1$ at $i$th position.
Let~$B_i^\prime$ denotes the $N\times{}(r-1)$-matrix~$B_i$ with last column removed.

Our construction is as follows.
Consider the direct product 
$$
\R^N\times\R^{r-1}\simeq\R^{N+r-1},
$$
with $\R^N$ being understood as the ``vertical'', and $\R^{r-1}$ as the ``horizontal'' subspace;
the coordinates on these spaces being denoted by 
$y=(y_1,\ldots,y_N)$ and $x=(x_1,\ldots,x_{r-1})$, respectively.
At every point $\bar b\in\R^N$, consider the affine $(r-1)$-dimensional subspace of $\R^{N+r-1}$
through~$\bar b$ defined by
$$
y=B(\bar b)^\prime{}x+\bar b.
$$
It follows from the maximality of the rank  that
any two of these affine spaces are skew,
and we obtain an affine Hopf fibration on $\R^{N+r-1}$. 

Next we give geometric arguments for the nonexistence.
Consider the first example after dimension three, that of $\R^4$. 

\begin{ex}
The space $\R^4$ has no affine Hopf fibrations.
Indeed, the only affine subspaces of~$\R^4$ that can be pairwise skew are
one-dimensional.
However, if $p=1$ then $n-p=3$ is odd, and so $\rho(n-p)=1$.
The condition of Theorem~\ref{TheThm} is not satisfied.
The same holds true for all
spaces of dimension $2^k$, that is, the space $\R^{2^n}$ admits no affine Hopf fibrations.
\end{ex}

\begin{figure}[hbtp]
\vskip-0.5cm
\includegraphics[width=8cm]{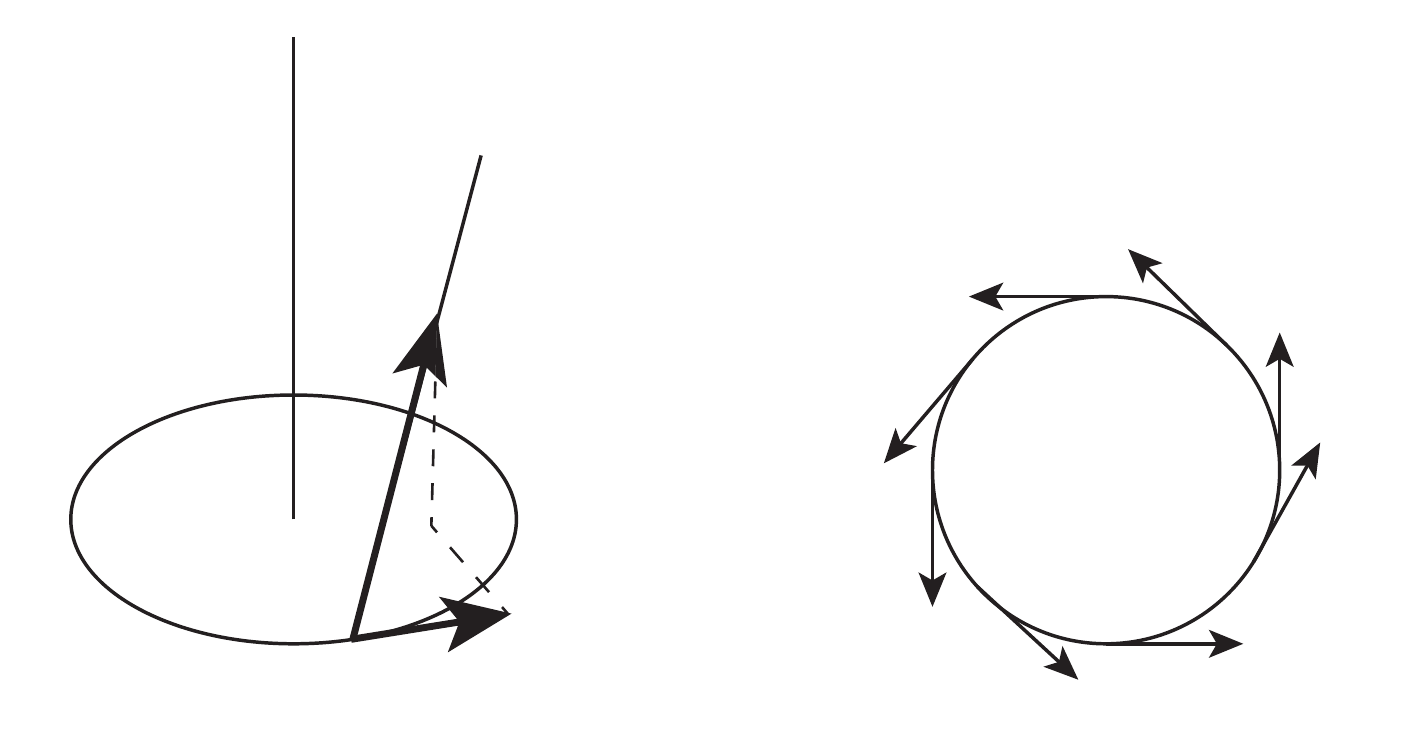}
\vskip-0.5cm
\caption{Affine Hopf fibrations and vector fields on spheres.}
\label{EzhFig}
\end{figure}

For the proof, consider first $\R^3$ equipped with the fibration depicted
in Figure~\ref{Hof}.
Choose any fiber, say vertical, and a two-dimensional plane orthogonal to it.
Draw a circle around the point of intersection of the plane and the fiber.
Pick any point $x$ of the circle and consider the unit  vector
that belongs to the fiber through~$x$.
Project this vector to the tangent line of the circle
(the projection goes in two steps: first project to the horizontal plane, then to the tangent),
see Figure~\ref{EzhFig}.
\noindent
Since the fiber through $x$ is skew with the
vertical line, the projection is non-zero.

Assume now that $\R^4$ is fibered by pairwise affine lines.
Exactly the same construction would then give a non-vanishing
vector field on a two-dimensional sphere that lies in a three-dimensional
affine subspace, orthogonal to a fiber.
This, however, contradicts the famous {\it hairy ball theorem}:  $S^2$ does not admit a non-vanishing continuous tangent vector field.
It follows that $\R^4$ admits no affine Hopf fibrations.

A  similar argument is used to deduce the ``necessary'' part of Theorem~\ref{TheThm} 
from Theorem \ref{Athm}. 

In fact, Theorem \ref{TheThm} extends to a semi-local result: 
its assertion holds for fibrations 
by pairwise skew affine subspaces of a neighborhood of a single fiber.

\begin{ex}
Let us now consider $p=8$ and $n=24$.
Since $\rho(16)=9$, the condition of Theorem~\ref{TheThm} is satisfied.
We therefore have an affine Hopf fibration  of $\R^{24}$ with 8-dimensional fibers.
This is the first example of
affine Hopf fibration which does not come from a projection of a classical Hopf fibration.
\end{ex}

More generally, a remarkable series of affine Hopf fibrations
corresponds to the following dimensions:
$$
\left(p,\,n-p\right)=
\left(\rho(2^n)-1,\,2^n\right).
$$
None of them is a projection of a classical Hopf fibration.

\subsection{A table}

The condition of Theorem~\ref{TheThm} is implicit (since $p$ appears in both parts of the
inequality).
Therefore, for a given $n$, the situation has to be analyzed.

Consider the cases of  $n\leq80$.
The dimensions $(p,n)$ for which affine Hopf fibrations of~$\R^n$ with 
$p$-dimensional fibers exist are presented in the following tables.

Recall the semi-local extension of Theorem \ref{TheThm}. It implies that, given a (local) $(p, n)$ affine Hopf fibration, the intersection with a hyperplane transverse to the fibers yields a (local) $(p - 1, n - 1)$ affine Hopf fibration. 
 This leads to the notion of a {\it dominant}  $(p, n)$ affine Hopf fibration: it is the case when a $(p + 1, n + 1)$  affine Hopf fibration does not exist. In the table below, the dominant pairs $(p, n)$ are shown in boldface.

$$
\setlength{\extrarowheight}{3pt}
\begin{array}{c||c|c|c|c|c|c|c|c|c|c|c|c|c|c|}
n  &{\bf 3} &4 &5 &6 &{\bf 7} &8 & 9&10 & 11 & 12&13&14&{\bf 15}
&16\\
\hline
\hline
p & {\bf 1} & & 1 & 2 &{\bf 3} &  &1 & 2 & 3      & 4 &5&6&{\bf 7} & \\
   &             & &    &   &{\bf 1}  &  &    &   &{\bf 1}&   & 1&2&{\bf 3}&    \\
   &&&&&&&&&&&&&{\bf 1}&
\end{array}
$$ 
The fibrations with $n=3,\,7,\,15$ are the central projections of the Hopf fibrations.
The next values are as follows.
$$
\setlength{\extrarowheight}{3pt}
\begin{array}{c||c|c|c|c|c|c|c|c|c|c|c|c|c|c|c|c|c}
n  &17&18&19&20&21&22&23&{\bf 24}&25 &26&27&28 &29&30 &{\bf 31}&32\\
\hline
\hline
p & 1&2&3       &4&5&6&7        &{\bf 8}&1 & 2 & 3   & 4 & 5& 6 & {\bf 7} &&\\
    &   &  &{\bf 1}&  &1&2&{\bf 3}&     &  &  & {\bf 1}&   & 1& 2 &{\bf 3}&  &  \\
   &&&&&&&{\bf 1}&&&&&&&&{\bf 1}&&
\end{array}
$$ 
including two interesting cases: $(p,n)=(8,24)$ and $(7,31)$.
The next values are:
$$
\setlength{\extrarowheight}{3pt}
\begin{array}{c||c|c|c|c|c|c|c|c|c|c|c|c|c|c|c|c|c}
n  &33& 34&35&36&37&38
&39&40&{\bf 41}&42&43&44&45 &46&{\bf 47}&48 \\
\hline
\hline
p&1&  2 &3       &4  &5&6 & 7        &8&{\bf 9}&2 & 3&4& 5&6 & {\bf 7 }& \\
     &&     &{\bf 1}& &1   & 2 &{\bf 3}&  &{1}&   &{\bf 1}&  &1&2& {\bf 3} &  \\
 &&&&&&&{\bf 1}&&&&&    &&   & {\bf 1}&
\end{array}
$$ 
We see another example of an interesting ``non-Hopf'' fibration: $\R^9\to\R^{41}\to\R^{32}$.
The next values are as follows.
$$
\setlength{\extrarowheight}{3pt}
\begin{array}{c||c|c|c|c|c|c|c|c|c|c|c|c|c|c|c|c|}
n &49&50 &51&52&53& 54&55&{\bf 56}&57&58
&59&60&61&62&{\bf 63}&64\\
\hline
\hline
p& 1&2& 3    &4 & 5&6&7        &{\bf 8}&1 &2&3       &4&5&6 &{\bf 7}&  \\
   & & &{\bf 1}&   &1 &2&{\bf 3}&          &   &   &{\bf 1}&  &1&2&{\bf 3}   &  \\
   &&&&&&&{\bf 1}&&&&&&&&{\bf 1}&
\end{array}
$$ 
and
$$
\setlength{\extrarowheight}{3pt}
\begin{array}{c||c|c|c|c|c|c|c|c|c|c|c|c|c|c|c|c|c|c|c|}
n&65 &66&67&68 &69&70 &71&72&73&74&{\bf 75}&76&77&78
&{\bf 79}&80\\
\hline
\hline
p&1&2& 3& 4& 5&6& 7      &8  &9&10&{\bf 11}&4&5&6 &{\bf 7}&  \\
 & &  & 1&   & 1 &2&{\bf 3}&   &1&2  &{3}  &    &1&2&{\bf 3}&    \\
 &  &   & &&&&{\bf 1}&&&&{\bf 1}&&&&{\bf 1}&
\end{array}
$$ 
We observe several more non-Hopf fibrations with  relatively large fibers.

Another observation is that some of the spaces
admit no fibrations with skew fibers at all.
In addition to the cases when the dimension is a power of 2, this happens for  $\R^{80}$.


\subsection{The complex case}

The definition of affine Hopf fibrations in $\C^n$ is exactly the same as in the real case.
Surprisingly, the problem is much harder
in this case, and seems not to be reducible to linear algebra.

Only necessary conditions are known for the existence of 
complex affine Hopf fibrations.
It was proved in \cite{Har} that an affine Hopf fibration of $\C^n$ with fibers of
dimension $p$ may exist only if, for each integer $r$ with
$0\leq{}r\leq{}p$, the coefficients of $t^r$ in the power series expansion of
$$
\left(
\frac{t}{\ln(1+t)}
\right)^{n-p}
$$
is an integer.

The minimal dimension $n$, considered a function of $p$, grows very  fast. 
For example, for a complex skew fibration with $1$-dimensional fibers,
 the ambient space is at least $3$-dimensional (like in the real case), 
 but for $2$-dimensional fibers, the ambient dimension is not less than $26$, and 
for $4$-dimensional fibers,  the ambient dimension is not less than $2884$.

\bigskip

{\bf Acknowledgments.}
We are grateful to
D. B.~Fuchs, M.~Harrison, and S.~Morier-Genoud
for enlightening discussions,
and to L.~Shenderova-Fock for a nice picture.
The second author was partially supported by the NSF grants DMS-1105442 and  DMS-1510055. 


\end{document}